\documentclass{amsart}
\usepackage{graphicx}
\usepackage{latexsym}
\usepackage{amsfonts}
\usepackage[all]{xy}
\usepackage{amssymb, mathrsfs, amsfonts, amsmath}
\usepackage{amsbsy}
\usepackage{amsfonts}
\newtheorem{corollary}{Corollary}

\newtheorem{theorem}{Theorem}
\newtheorem{example}{Example}
\numberwithin{equation}{section}

\begin{document}
	\title[An analysis of symmetry groups of generalized $m$-quasi-Einstein manifolds]{An analysis of symmetry groups of generalized $m$-quasi-Einstein manifolds}
	
	\author{Paula Correia$^{1}$}
	
	\author{Benedito Leandro$^{2}$}
	
	\author{Romildo Pina$^{3}$}

	\address{$^{1,\,2,\,3}$ IME, Universidade Federal de Goi\'as,
		Caixa Postal 131, CEP 74001-970, Goi\^ania, GO, Brazil.}

	\email{paulacorreiacatu@hotmail.com$^{1}$}
	\email{bleandroneto@ufg.br$^{2}$}
	\email{romildo@ufg.br$^{3}$}

	\date{\today}

	\begin{abstract}
		\noindent
		In this paper emphasis is placed on how the behavior of the
		solutions of a PDE is affected by the geometry of the generalized $m$-quasi-Einstein manifold, and vice
		versa. Considering a $n$-dimensional generalized $m$-quasi-Einstein manifold which is conformal to a pseudo-Euclidean space, we prove the most general symmetry group of maximal dimension. Moreover, we demonstrate that there is no different low dimensional invariant on a generalized $m$-quasi-Einstein manifold. As an application, we use the invariant structure of the metric to provide an example of shrinking $m$-quasi-Einstein manifold (cf. Example 3). A discussion about the fluid ball conjecture was made.
		\end{abstract}
	
	\maketitle
	\section{Introduction}

	An $n$-dimensional Riemannian manifold $(M^n,g)$ is generalized $m$-quasi-Einstein if there exist two smooth functions $f$ and $\lambda$ on $M$ such that
	\begin{equation}\label{eqp1-1}
	Ric_g + Hess_gf - \frac{1}{m} df\otimes df = \lambda g,
	\end{equation}
	where $m\in(0,+\infty]$. When $m\in(0,+\infty)$, we can make the change $ h = e^{-\frac{f}{m}} $ and get the equation
	\begin{equation} \label{eq1}
	Ric_g - \frac{m}{h}Hess_gh = \lambda g.
	\end{equation}

	In \cite{Catino}, Catino introduced the notion of generalized quasi-Einstein manifolds. He proved that a complete generalized quasi-Einstein manifold with
	harmonic Weyl tensor and vanishing radial Weyl curvature is locally a warped product with $(n-1)-$dimensional Einstein fiber.

	The importance of understanding and giving explicit solutions for generalized $m$-quasi-Einstein manifolds arises from the fact that they are closely related to Einstein warped product manifolds (cf. \cite{Besse}). Furthermore, it is well known that \eqref{eqp1-1} generalizes the notion of gradient Ricci solitons, and Einstein manifolds. Also, \eqref{eq1} generalizes several important metrics, e.g., critical metrics and static metrics (cf. \cite{Leandro} and the references therein). Therefore, this problem has great importance in physics.
	
	In this paper we give a simple proof for the fluid ball conjecture (cf. \cite{masood}) considering that the spatial factor of the perfect fluid is conformally flat. This conjecture was proposed by Yau in the 1982 list of unsolved problems in General Relativity (cf. \cite{yau}). The fluid ball conjecture states that a non-rotating stellar model is spherically symmetric. Here, we prove that a perfect fluid space-time in which the spatial factor is conformal to a pseudo-Euclidean space is either invariant by the action of a pseudo-orthogonal group or by the action of a translation group.
	
	Throughout history, several methods of reduction (ansatz) of PDEs were used in differential geometry to provide examples or even full classification of metrics (cf. \cite{barbosaketipina,MarceloBeneditoRomildo,BeneditoeJoaoPaulo,MarcioeRomildo}). As examples of ansatz methods we refer to the Lie group theory and the method of characteristics (cf. \cite{olver1,olver2}). The reduction method used in this paper is based on those two (cf. Theorem \ref{teoremaprincipal} and Theorem \ref{teoremadenaoexistencia}). It is worth to point out that the uniqueness of static vacuum Einstein space-time was provide using this tecnique (cf.  \cite{joao}). 
	
	Recently, Ribeiro Jr and Tenenblat \cite{RibeiroTenenblat} classified the $n$-dimensional $m$-quasi-Einstein manifolds invariants under the action of a translation group, in which $\lambda$ is constant. They also gave a complete classification when $\lambda=0$, $m\geq1$ or $m=2-n$. 
	
	In this work we study conformally flat generalized $m$-quasi-Einstein manifolds satisfying \eqref{eq1}. However, unlike \cite{RibeiroTenenblat} we do not fix any kind of symmetry. In fact, we completely describe the most general ansatz capable of reducing the system of PDEs obtained from \eqref{eq1} to a system of ODEs (ordinary differential equations).  Moreover, as far as we know there is no classification for conformally flat pseudo-Euclidean generalized $m$-quasi-Einstein manifolds. 

	Hereafter, we will establish the needed notations to announce our main results. Let $(\mathbb{R}^{n}, g)$ be the standard pseudo-Euclidean space with coordinates $(x_{1}, \cdots, x_{n})$ and metric components $g_{ij} = \delta_{ij}\varepsilon_{i}$, $1\leq i, j\leq n$, where 
	$\varepsilon_{i} = \pm1$. We want to find smooth functions $\varphi$, $h$ and $\lambda$ defined on an open subset $\Omega \subset \mathbb{R}^n$ such that, for $\bar{g}$ given by
	\begin{eqnarray*}\label{gbar}
		\bar{g}=\frac{g}{\varphi^2},
	\end{eqnarray*}
	$(\Omega, \bar{g})$ is a generalized $m$-quasi-Einstein manifold with potential function $h$, i.e.,
	\begin{eqnarray}
	Ric_{\bar{g}}-\dfrac{m}{h}Hess_{\bar{g}}(h)=\lambda\bar{g}, \label{mqeg-bar}
	\end{eqnarray}
	where $m\in(0,+\infty)$, $Ric_{\bar{g}}$ and $Hess_{\bar{g}}(h)$ are, respectively, the Ricci tensor and the Hessian of the metric $\bar{g}$.
	
	In general, find explicit solutions for the system of partial differential equations (PDE) that describes this problem in local coordinates is not simple. However, by an ansatz method we are able to transform the PDE system into a system of ordinary differential equations (ODE). Under this approach we can find explicit solutions the generalized $m$-quasi-Einstein equation (see \cite{olver1,olver2}). Moreover, we will prove the most general form that a smooth function $\xi:\Omega\subseteq\mathbb{R}^{n}\rightarrow\mathbb{R}$ such that $h\circ\xi$, $\varphi\circ\xi$ and $\lambda\circ\xi$ satisfying \eqref{mqeg-bar} must have. 
	
	Our first result tells us that a generalized $m$-quasi-Einstein manifold reduced by our ansatz method is invariant under the action of tô the pseudo-orthogonal group and the action of the translation group. Moreover, Theorem \ref{teoremadenaoexistencia} demonstrates that there is no other group-invariant even in low dimensions.
	
	It is worth pointing out that in \cite{TenenblatWinternitz} the intrinsic generalized wave and sine-Gordon equations were studied assuming the same ansatz method. The authors provided the symmetry groups of local point transformations for those equations, and then they used those invatiant groups to obtain classes of exact solutions.
	
	Also in \cite{joaoromildo} the authors used the Lie-point symmetries to find metrics that would solve the Ricci curvature and the Einstein equations. They provided a large class of group-invariant solutions and examples of complete metrics defined globally in $\mathbb{R}^{n}$.
	
	Here, using a different approach, we will describe all the invariant groups for a generalized $m$-quasi-Einstein manifold conformal to a pseudo-Euclidean space, based upon our ansatz method. Furthermore, it is important to say that Theorem \ref{teoremaprincipal} and Theorem \ref{teoremadenaoexistencia} also can be applied to $m$-quasi-Einstein manifolds.
	
	Without further ado, we state our main results.

\begin{theorem} \label{teoremaprincipal}
	Let $(\mathbb{R}^{n}, g)$ be the standard pseudo-Euclidean space, $n\geq 2$, with Cartesian coordinates $x=(x_{1}, ..., x_{n})$, $g_{ij}=\delta_{ij}\varepsilon_i$ and let $\Omega\subseteq \mathbb{R}^n$ be an open subset. Consider non-constant smooth functions $h(\xi),\,\varphi(\xi),\, \lambda(\xi): \Omega\subseteq\mathbb{R}^{n} \rightarrow \mathbb{R}$. Then
	\begin{eqnarray*}
	\displaystyle \left(\mathbb{R}^{n},\,\bar{g}=\frac{1}{\varphi^{2}(\xi)}g,\,h(\xi),\,\lambda(\xi)\right)
	\end{eqnarray*}
	is a generalized $m$-quasi-Einstein manifold satisfying \eqref{mqeg-bar} if, and only if, the function $\xi : \Omega\subseteq\mathbb{R}^{n} \rightarrow \mathbb{R}$ is of the form
	\begin{eqnarray}\label{invar}
	\xi=P\left(\sum_{k=1}^{n}a\varepsilon_kx_k^2+b_kx_k+c_k\right),
	\end{eqnarray}
	where $a, b_k, c_k\in \mathbb{R}$ and $P$ is at least a $C^1$ function.
\end{theorem}
	
	Now we show a result concerning the rigidity of generalized $m$-quasi-Einstein metrics. We prove that there is no other symmetry group of low dimension on a generalized $m$-quasi-Einstein manifold.
	This question was mentioned also in \cite{BeneditoeJoaoPaulo}, where the authors proved all the maximal invariant groups for the gradient Ricci soliton by the same ansatz method used in this work. However, they did not provide the low dimensional symmetries for the gradient Ricci solitons.
	
	The nonexistence of $m$-quasi-Einstein metrics is an important issue. As we have said before, it can indicate the impossibility to construct Einstein warped products metrics or even gradient Ricci solitons. On this subject, we recommend the reader to see \cite{case}.

	\begin{theorem} \label{teoremadenaoexistencia}
		Let $(\mathbb{R}^{n}, g)$ be the standard pseudo-Euclidean space, $n\geq 2$, with Cartesian coordinates $x=(x_{1}, ..., x_{n})$, $g_{ij}=\delta_{ij}\varepsilon_i$.
		Considering smooth functions 
		$h(\xi),\,\varphi(\xi)$ and $\lambda(\xi)$ satisfying \eqref{mqeg-bar} where $$\displaystyle \xi=\xi(x_{1},\,\ldots,\,x_{n-1}).$$
		Then,  $$\xi=P\left(\displaystyle\sum_{k=1}^{n-1}a_kx_{k}+b_{k}\right),$$
		where $a_k,\,b_{k}\in\mathbb{R}$ and $P$ is at least a $C^1$ function.
		Moreover, if the invariant has a different form then the generalized $m$-quasi-Einstein manifold is trivial, i.e., either $\varphi$ or $h$ is a constant functions.
	\end{theorem}

Consequently we have the following classification theorem for a generalized $m$-quasi-Einstein manifold.

\begin{theorem}
	Consider a conformally flat semi-Riemannian generalized $m$-quasi-Einstein manifold invariant under the action of some symmetry group. Then this group is either a pseudo-orthogonal group or a translation group.
	\end{theorem}
	
	For instance, we can gather that $\mathbb{R}\times\mathbb{S}^{n-1}$ do not represent a group-invariant for a generalized $m$-quasi-Einstein manifold. Precisely, we obtain the next result.
	\begin{corollary} \label{coro1}
		Let $(\mathbb{R}^{n}, g)$ be the standard pseudo-Euclidean space, $n\geq 2$, with Cartesian coordinates $x=(x_{1}, ..., x_{n})$, $g_{ij}=\delta_{ij}\varepsilon_i$.
		Considering $h(\xi),\,\varphi(\xi)$ and $\lambda(\xi)$ smooth functions satisfying \eqref{mqeg-bar} where $$\displaystyle \xi=
		\sum_{k=1}^{n-1}\varepsilon_kx_k^2.$$
		Then, 	$$\displaystyle \left(\mathbb{R}^{n},\,\bar{g}=\frac{1}{(\varphi^{2}(\xi)}g,\,h(\xi),\,\lambda(\xi)\right)$$ is a trivial generalized $m$-quasi-Einstein manifold, i.e., either $\varphi$ or $h$ is a constant functions.
	\end{corollary}

Let us give a physical application of our main results, Theorem \ref{teoremaprincipal} and Theorem \ref{teoremadenaoexistencia}.

\begin{corollary}
	Considering a $n$-dimensional Riemannian manifold that is a generalized $m$-quasi-Einstein satisfying \eqref{eq1} such that $m=1$, $\lambda = \frac{(\mu-\rho)}{(n-1)}$ and 
	\begin{eqnarray*}\label{energypressure}
	\Delta h=\left(\dfrac{n-2}{(n-1)}\mu+\dfrac{n}{n-1}\rho\right)h,
	\end{eqnarray*}
	where $\mu,\,\rho:M\rightarrow\mathbb{R}$ are smooth functions called density and pressure, respectively, we have a static perfect fluid equation. Under the conditions of Theorem \ref{teoremaprincipal} and Theorem \ref{teoremadenaoexistencia} the spatial factor of the perfect fluid space-time is either invariant by the action of a pseudo-orthogonal group or by the action of a translation group. 
\end{corollary}
	
	In what follows, we provide the reduction of the PDE system \eqref{mqeg-bar} into a ODE, which is a consequence of Theorem \ref{teoremaprincipal}.
	
	\begin{theorem} \label{corolario}
		Under the same conditions of Theorem \ref{teoremaprincipal}, for any function $\varphi(\xi)$, $$\displaystyle \left(\mathbb{R}^{n},\bar{g}=\frac{1}{\varphi^2}g\right)$$ is a generalized $m$-quasi-Einstein manifold if, and only if, the function $h$ is the solution of the ordinary differential equation
		\begin{equation} \label{equacaodeh}
		(n-2)h\varphi''-m\varphi h''-2m\varphi'h'=0,
		\end{equation}
		where the function $\lambda$ is given by
		\begin{eqnarray} \label{equacaodelambda}
		\lambda &=& 2a\varphi\left[(n-2)\varphi'-m\varphi \frac{h'}{h}\right]\nonumber\\
		&+&\left[\varphi\varphi''-(n-1)(\varphi')^2+m\varphi\varphi'\frac{h'}{h}\right](4a\xi+S)+2na\varphi\varphi'.
		\end{eqnarray}
		Here  $S=\sum_{k=1}^{n}(\varepsilon_kb_k-4ac_k)$ and $a,\,b_k,\,c_k\in\mathbb{R}$ for all $1\leq k\leq n$.
	\end{theorem}

	Now we explicitly provide families of solutions for generalized $m$-quasi-Einstein manifolds. The first one is invariant by rotations and the second example is invariant by translations. The last example is invariant by rotations and shows that, for a given choice of the conformal factor, the ODE system in Theorem \ref{corolario} has explicit solution.
	
	\begin{example}
		Let us consider $r=\sum_{k=1}^{n}x_k^2$ and the function $\varphi(r) = e^{\alpha r+\beta}$, with $\alpha,\,\beta\in\mathbb{R}$. From Theorem \ref{corolario} we have that $h$ is the solution of 
		$$h'' + 2\alpha h' - \frac{(n-2)}{m}\alpha^2h=0.$$
		Then, if $m>n-2$,
		$$h(r) = c_1e^{r_1r}+c_2e^{r_2r},$$
		where $r_1=-\alpha+|\alpha|\sqrt{\frac{m-(n-2)}{m}}$, $r_2=-\alpha-|\alpha|\sqrt{\frac{m-(n-2)}{m}}$ e $c_1,c_2\in\mathbb{R}$. Moreover,
		$$\lambda(r) = e^{2(\alpha r+\beta)}\left[-4(n-2)\alpha^2 r+4(n-1)\alpha+m(4\alpha r-2)\frac{c_1r_1e^{r_1 r}+c_2r_2e^{r_2 r}}{c_1e^{r_1r}+c_2e^{r_2 r}}\right].$$
		In this case, the solutions are globally defined.  Taking $\varphi(r) = e^{-r}$, we have that $\varphi$ is bounded and the metric $\bar{g}$ is complete.
	\end{example}
	
	\
	
	\begin{example}
		Considering $\xi=\sum_{k=1}^{n}b_kx_k$ and the function $\varphi(\xi) = e^{a\xi+b}$, with $a,b\in\mathbb{R}$. By Theorem \ref{corolario} we have that if $m>n-2$, $h$ is given by
		$$h(\xi) = c_1e^{r_1\xi}+c_2e^{r_2\xi},$$
		where $r_1=-a+|a|\sqrt{\frac{m-(n-2)}{m}}$, $r_2=-a-|a|\sqrt{\frac{m-(n-2)}{m}}$ e $c_1,c_2\in\mathbb{R}$. Therefore,
		$$\lambda(\xi) = \varepsilon_{i_{0}}ae^{2(a\xi+b)}\left[m\frac{c_1r_1e^{r_1\xi}+c_2r_2e^{r_2\xi}}{c_1e^{r_1\xi}+c_2e^{r_2\xi}}-(n-2)a\right],$$
		where $\varepsilon _{i_ 0}=\displaystyle\sum_{k}\varepsilon_{k}b_{k}^{2}$. In this case, the solutions are globally defined. We remark that if $b=\displaystyle\sum_{k}b_{k}\frac{\partial}{\partial x_{k}}$ is a null vector (lightlike) we have $\varepsilon _{i_ 0}=0$, and therefore $\lambda=0$, i.e., a steady quasi-Einstein manifold.
	\end{example}

	\

	\begin{example}
		
		Consider $\displaystyle r=\sum_{k=1}^{n}x_k^2$, and $\varphi=\sqrt{r}$, such that $r>0$, in Theorem \ref{corolario}. By \eqref{equacaodeh}, $h$ is a solution of the Euler equation
		$$r^2h''+rh'+\frac{(n-2)}{4m}h=0.$$
		Since $n\geq 2$, $h$ is given by
		\begin{equation*}
		h(r)=
		\begin{cases}
		c_1+c_2\log r, & n=2,\\
		c_1\sin(\mu\log r)+c_2\cos(\mu\log r), & n>2,
		\end{cases}
		\end{equation*}
		with $c_{i}\in\mathbb{R}$ and $\displaystyle \mu=\frac{1}{2}\sqrt{\frac{n-2}{m}}$. Note that, by equation \eqref{eq1}, we should have $h(r) \neq 0$ for all $r$. This function must be totally positive (or negative), and for $n\geq 2$, $h$ is not globally defined.
		
		Therefore, if $n=2$ we have a steady $m$-quasi-Einstein manifold, and if $n>2$ we have a shrinking $m$-quasi-Einstein manifold in which $\lambda=(n-2)$, both non complete.
	\end{example}
	
	\vspace{.2in}

	\section{Background}
	\vspace{.2in}
	We denote $\varphi,_i$ and $h,_i$ the first order derivatives, and $\varphi,_{ij}$ and $h,_{ij}$ as the second order derivatives of the functions $\varphi$ and $h$ with respect to $x_{i}$ and $x_{i}x_{j}$, respectively.
	
	Let $(\mathbb{R}^{n}, g)$ is the pseudo-Euclidean space with coordinates $x=(x_{1}, ..., x_{n}),  \ g_{ij}=\delta_{ij}\varepsilon_i$. From the conformal structure (see \cite{Besse}), if $\displaystyle \bar{g}=\frac{1}{\varphi^{2}}g $, we obtain
	\begin{equation} \label{eq9}
	Ric_{\bar{g}}- Ric_{g}=\frac{1}{\varphi^{2}}\left\lbrace \left( n-2\right)\varphi Hess_{g}\varphi + \left[ \varphi\Delta_{g}\varphi-\left( n-1\right)\left|\nabla_{g}\varphi \right|^{2}  \right]g\right\rbrace. 
	\end{equation}	
	
	Thus, for a tangent base $X_{1},...,X_{n}$ of $\mathbb{R}^{n}$ we get
	
	\begin{equation} \label{eq12}
	\begin{cases}
	\displaystyle
	\left(Hess_{\bar{g}}(h) \right)_{ij}=h,_{ij}+\frac{ \varphi,_{j}}{\varphi}h,_{i}+\frac{\varphi,_{i}}{\varphi}h,_{j}, \ \ i\neq j\\ \\
	\displaystyle
	\left(Hess_{\bar{g}}(h) \right)_{ii}=h,_{ii}+2\frac{\varphi,_{i}}{\varphi}h,_{i}-\varepsilon_i\sum_{k=1}^{n}\varepsilon_k\frac{\varphi,_{k}}{\varphi}h,_{k}, \ \ i=j,
	\end{cases}
	\end{equation}
	where  $Hess_{\bar{g}}(h)\left( X_{i}, X_{j}\right) = \left(Hess_{\bar{g}}(h) \right)_{ij}$.
	
	Remember that
	\begin{equation} \label{eq10}
	Ric_{\bar{g}}-\frac{m}{h} Hess_{\bar{g}}h=\lambda\bar{g}, \ \ \ \lambda \in C^{\infty} (\mathbb{R}^{n}).
	\end{equation} 
	Replacing (\ref{eq9}) and (\ref{eq12}) in (\ref{eq10}), provided that $\displaystyle \Delta_g\varphi=\sum_{k=1}^{n}\varepsilon_k\varphi,_{kk}$ and $\displaystyle |\nabla_g\varphi|^2=\sum_{k=1}^{n}\varepsilon_k\varphi,_k^2$, we get 
	
	\begin{eqnarray} \label{EDP}
	\left\{
	\begin{array}{lcc}
	(n-2)h\varphi,_{ij}-m\left(\varphi h,_{ij}+\varphi,_jh,_i+\varphi,_ih,_j \right)=0; \ \forall\, i\neq j,\\\\
	\displaystyle (n-2)h\varphi\varphi,_{ii} + \varepsilon_ih\sum_{k=1}^{n}\varepsilon_k\left[\varphi\varphi,_{kk}-(n-1)\varphi,_k^2\right]\\\\
	\quad\quad-m\left[\varphi^2h,_{ii}+2\varphi\varphi,_ih,_i-\varepsilon_i\sum_{k=1}^n\varepsilon_k\varphi\varphi,_kh,_k\right]=\varepsilon_i\lambda h; \ \mbox{for}\quad i=j.
	\end{array}\right.
	\end{eqnarray}
	
	\vspace{.2in}

	\section{Proof of the Main Results}

	\vspace{.2in}

	\noindent\textbf{Proof of Theorem \ref{teoremaprincipal}:}
	Consider that \eqref{EDP} admits non-trivial solutions such that $f\circ\xi,\, h\circ\xi,\, \varphi\circ\xi$ and $\lambda\circ\xi$, where $\xi : \Omega\subseteq\mathbb{R}^n\rightarrow\mathbb{R}$ is a smooth function. Then, from the first equation of (\ref{EDP}) we obtain that
	
	\begin{equation} \label{eq1cominvariancia}
	[(n-2)h\varphi''-m\varphi h''-2m\varphi'h']\xi,_i\xi,_j + [(n-2)h\varphi'-m\varphi h']\xi,_{ij}=0.
	\end{equation}
	
	Note that $(n-2)h\varphi'-m\varphi h'\neq 0$. In fact, if $(n-2)h\varphi'-m\varphi h'=0$ then
	
	$$(n-2)\frac{\varphi'}{\varphi}= m\frac{h'}{h}.$$
	
	From the above identity we get
	\begin{equation}\label{igualdadedoMarcio}
	\frac{\varphi''}{\varphi}=\left(\frac{\varphi'}{\varphi}\right)'+\left(\frac{\varphi'}{\varphi}\right)^2,
	\end{equation}
	from \eqref{eq1cominvariancia} we can infer that either $\varphi$ or $h$ is constant, which is a contradiction. Thus, from \eqref{eq1cominvariancia} we get
	
	\begin{equation}\label{2.8}
	\frac{\xi,_{ij}}{\xi,_i\xi,_j}=\frac{-[(n-2)h\varphi''-m\varphi h''-2m\varphi'h']}{(n-2)h\varphi'-m\varphi h'}=F(\xi).
	\end{equation}
	
	So,
	$$\log(\xi,_i)=\int F(\xi)d\xi + F_i(\hat{x_j}),$$
	where the symbol $\hat{x_j}$ denotes that the function does not depend on $x_j$. That is,
	$$\xi,_i=e^{\int F(\xi)d\xi}e^{F_i(\hat{x_j})}.$$
	Since this is true for all $i\neq j$, denote by $G(\xi)=e^{\int F(\xi)d\xi}$ and $G_i(x_i)=e^{F_i(\hat{x_j})}$, therefore,
	\begin{equation} \label{2.14}
	\xi,_i=G(\xi)G_i(x_i).
	\end{equation}
	
	Now, contracting \eqref{eq1} we obtain
	\begin{equation*}
	(n-1)\sum_{k=1}^{n}\varepsilon_k\left(2\varphi\varphi,_{kk}-n\varphi,_k^2\right)-\frac{m}{h}\sum_{k=1}^{n}\varepsilon_k\left[\varphi^2h,_{kk}-(n-2)\varphi\varphi,_kh,_k\right]=n\lambda.
	\end{equation*}
	Multiplying this equation by $h$ and the second equation of \eqref{EDP} by $\varepsilon_{i}n$, we can conclude that
	
	$$[(n-2)h\varphi''-m\varphi h''-2m\varphi'h']\left(\varepsilon_in\xi,_i^2-\sum_{k=1}^{n}\varepsilon_k\xi,_k^2\right)$$
	$$+ [(n-2)h\varphi'-m\varphi h']\left(\varepsilon_in\xi,_{ii}-\sum_{k=1}^{n}\varepsilon_k\xi,_{kk}\right)=0.$$
	
	Thus, from \eqref{2.8} we have
	$$\varepsilon_in[\xi,_{ii}-F(\xi)\xi,_i^2]=\sum_{k=1}^{n}\varepsilon_k[\xi,_{kk}-F(\xi)\xi,_k^2].$$
	Hence,
	$$\varepsilon_in\left(\xi,_ie^{-\int Fd\xi}\right),_i=\sum_{k=1}^{n}\varepsilon_k\left(\xi,_ke^{-\int Fd\xi}\right),_k.$$
	
	Then, by \eqref{2.14} we can infer that
	$$\varepsilon_inG_i'=\sum_{k=1}^{n}\varepsilon_kG_k'.$$
	Since the left-hand side depends only on $x_i$, we have
	$$\varepsilon_iG_i'=\varepsilon_jG_j', \ \forall i\neq j.$$
	
	Thus $G_i(x_i)=2a\varepsilon_ix_i+b_i$, with $a,b_i\in\mathbb{R}$. Therefore, $\displaystyle \frac{\xi,_i}{G_i}=G$ implies that
	$$\frac{\xi,_i}{2a\varepsilon_ix_i+b_i}=\frac{\xi,_j}{2a\varepsilon_jx_j+b_j},$$
	where we conclude that $\xi$ is of the form \eqref{invar}.
	
	The reverse statement is a straightforward computation.

	\hfill $\Box$

	\vspace{.2in}
	
	\noindent\textbf{Proof of Theorem \ref{teoremadenaoexistencia}:}
	Now, considering $\xi=\xi(x_1,\cdots,x_{n-1})$, for all $1\leq i,j \leq n-1$ we have
	
	\begin{center}
		\begin{tabular}{ccc}
			$\varphi,_i=\xi,_i\varphi'$ & $\varphi,_{ij}=\xi,_i\xi,_j\varphi''+\xi,_{ij}\varphi'$ & $\varphi,_{ii}=\xi,_i^2\varphi''+\xi,_{ii}\varphi'$ \\
			$h,_i=\xi,_ih'$ & $h,_{ij}=\xi,_i\xi,_jh''+\xi,_{ij}h'$ & $h,_{ii}=\xi,_i^2h''+\xi,_{ii}h'$ \\
		\end{tabular}
	\end{center}
	and  $$\varphi,_n=\varphi,_{in}=\varphi,_{ni}=h,_n=h,_{in}=h,_{ni}=0,\quad\mbox{for all}\quad i=1,\cdots,n.$$
	
	Then, from the first equation of \eqref{EDP}, for $i\neq j \neq n$, we obtain
	\begin{equation}\label{eqparaidiferentedej}
	[(n-2)h\varphi''-m\varphi h''-2m\varphi'h']\xi,_i\xi,_j + [(n-2)h\varphi'-m\varphi h']\xi,_{ij}=0.
	\end{equation}
	Considering $i=n$ or $j=n$, the first equation of \eqref{EDP} is trivially satisfied.
	
	Now, from the second equation of \eqref{EDP}, for $i\neq n$ we obtain
	\begin{eqnarray} \label{segundaeqdiferenteden}
	&&\varphi[(n-2)h\varphi''-m\varphi h''-2m\varphi'h']\xi,_i^2 + \varphi[(n-2)h\varphi'-m\varphi h']\xi,_{ii}\nonumber\\
	&+& \varepsilon_i\sum_{k=1}^{n-1}\varepsilon_k\{[h\varphi\varphi''-(n-1)h(\varphi')^2+m\varphi\varphi'h']\xi,_k^2 + h\varphi\varphi'\xi,_{kk}\}=\varepsilon_i\lambda h.
	\end{eqnarray}
	On the other hand, for $i=n$ we get
	\begin{equation*}
	\sum_{k=1}^{n-1}\varepsilon_k\{[h\varphi\varphi''-(n-1)h(\varphi')^2+m\varphi\varphi'h']\xi,_k^2 + h\varphi\varphi'\xi,_{kk}\}=\lambda h.
	\end{equation*}
	Using the above identity in  \eqref{segundaeqdiferenteden} leads us to
	\begin{equation}\label{eqparai=j}
	[(n-2)h\varphi''-m\varphi h''-2m\varphi'h']\xi,_i^2 + [(n-2)h\varphi'-m\varphi h']\xi,_{ii}=0.
	\end{equation}
	
	From \eqref{eqparaidiferentedej} and \eqref{eqparai=j} follows that if $(n-2)h\varphi'-m\varphi h'\neq 0$,
	\begin{equation} \label{igualdadeporF}
	\frac{\xi,_{ij}}{\xi,_i\xi,_j}=\frac{-[(n-2)h\varphi''-m\varphi h''-2m\varphi'h']}{(n-2)h\varphi'-m\varphi h'}=\frac{\xi,_{ii}}{\xi,_i^2}.
	\end{equation}
	Denoting by
	$$F(\xi)=\frac{-[(n-2)h\varphi''-m\varphi h''-2m\varphi'h']}{(n-2)h\varphi'-m\varphi h'},$$
	from the first equality of \eqref{igualdadeporF} we have
	$$\log\xi,_i=\int F(\xi)d\xi+F_i(\hat{x_j}), \ \forall i\neq j.$$
	Hence,
	\begin{equation} \label{daprimeiraigualdade}
	\xi,_i=e^{\int F(\xi)d\xi}e^{F_i(\hat{x_j})}=G(\xi)G_i(x_i).
	\end{equation}
	Similarly, by the second equality of \eqref{igualdadeporF},
	$$\log\xi,_i=\int F(\xi)d\xi+L_i(\hat{x_i}), \ \forall i.$$
	Then,
	\begin{equation} \label{dasegundaigualdade}
	\xi,_i=e^{\int F(\xi)d\xi}e^{L_i(\hat{x_i})}=G(\xi)K_i(\hat{x_i}).
	\end{equation}
	From \eqref{daprimeiraigualdade} and \eqref{dasegundaigualdade} we concluded that $G_i(x_i)=K_i(\hat{x_i})=a_i$, $a_i$ constant.
	
	Thus, we get
	$$\frac{\xi,_i}{a_i}=\frac{\xi,_j}{a_j}.$$
	Then, for all $i\neq j \neq n$, the characteristic for the above equation implies that
	$$\xi=P\left(\sum_{k=1}^{n-1}a_kx_k+b_k\right),$$
	where $a_k,\,b_k\in\mathbb{R}$ and $P$ is at least a $C^1$ function.
	
	Now, if $(n-2)h\varphi'-m\varphi h'=0$, then
	$$(n-2)\frac{\varphi'}{\varphi}=m\frac{h'}{h}.$$
	Using (\ref{igualdadedoMarcio}), from \eqref{eqparaidiferentedej} or \eqref{eqparai=j} we conclude that either $\varphi$ or $h$ is constant.

	\hfill $\Box$

	\vspace{.2in}
	
	\noindent\textbf{Proof of Theorem \ref{corolario}:}
	From Theorem \ref{teoremaprincipal}, we can suppose that
	$$\xi = \sum_{k=1}^{n}U_k(x_k).$$
	Therefore,
	$$\varphi,_i= \varphi'U_i', \hspace{1cm} \varphi,_{ij}=\varphi''U_i'U_j', \hspace{1cm} \varphi,_{ii}=\varphi''(U_i')^2+\varphi'U_i'',$$
	$$h,_i= h'U_i', \hspace{1cm} h,_{ij}=h''U_i'U_j', \hspace{1cm} h,_{ii}=h''(U_i')^2+h'U_i''.$$
	
	Thus, from the first equation of (\ref{EDP}) we obtain
	$$[(n-2)h\varphi''-m(\varphi h'' +2\varphi'h')]U_i'U_j'=0,$$
	from which we conclude that
	\begin{equation} \label{eq1cominvariante}
	(n-2)h\varphi''-m(\varphi h'' +2\varphi'h')=0.
	\end{equation}
	
	From the second equation of (\ref{EDP}) we have
	\begin{eqnarray*}
		&&\varphi[(n-2)h\varphi''-m(\varphi h'' +2\varphi'h')](U_i')^2+\varphi[(n-2)h\varphi'-m\varphi h']U_i''\nonumber\\
		&+&\varepsilon_i\sum_{k=1}^{n}\varepsilon_k\{[h\varphi\varphi''-(n-1)h(\varphi')^2+m\varphi\varphi'h'](U_k')^2+h\varphi\varphi'U_k''\}=\varepsilon_i\lambda h.
	\end{eqnarray*}
	Note that $\displaystyle \sum_{k=1}^{n}\varepsilon_k(U_k')^2=4a\xi+S$ and $\sum_{k=1}^{n}\varepsilon_kU_k''=2na$, where $S=\sum_{k=1}^{n}(\varepsilon_kb_k-4ac_k)$. Hence, using the equation \eqref{eq1cominvariante} we obtain
	\begin{eqnarray}\label{eq2cominvariante}
	&&2a\varphi[(n-2)h\varphi'-m\varphi h']\nonumber\\
	&+&[h\varphi\varphi''-(n-1)h(\varphi')^2+m\varphi\varphi'h'](4a\xi+S)+2nah\varphi\varphi'=\lambda h.
	\end{eqnarray}
	
	Finally, from  (\ref{eq1cominvariante}) and (\ref{eq2cominvariante}) the result follows.
	\hfill $\Box$
	
	\vspace{.2in}


	\
	
	\bibliographystyle{acm}
	\bibliography{bibliography}
	
\end{document}